\documentclass[12pt,a4paper]{amsart}
\usepackage[english]{babel}
\usepackage[utf8]{inputenc}
\usepackage{fancyhdr}
\usepackage{xcolor}

\usepackage{amssymb}
\usepackage{amsmath} 
\usepackage{amsthm} 
\usepackage{enumitem}
\usepackage{bbm}
\usepackage{bbold}

\newtheorem{thm}{Theorem}
\newtheorem{conj}{Conjecture}
\newtheorem{prop}{Proposition}

\newtheorem{cor}{Corollary}

\theoremstyle{proposition}
\newtheorem*{prop*}{Proposition}
\theoremstyle{proposition}
\newtheorem*{rem*}{Remark}

\numberwithin{lemma}{section}
\numberwithin{equation}{section}
\numberwithin{prop}{section}
\numberwithin{re}{section}

\title{Weighted value distributions of the Riemann zeta function on the critical line}
\author{Alessandro Fazzari}
\address{Universit\`a  di Genova,  Dipartimento di Matematica. 
 Via Dodecaneso 35, 16146 Genova, Italy}
\email{fazzari@dima.unige.it}

\DeclareMathOperator{\meas}{meas}

\begin{document}
\maketitle

\begin{abstract} We prove a central limit theorem for $\log|\zeta(1/2+it)|$ with respect to the measure $|\zeta^{(m)}(1/2+it)|^{2k}dt$ ($k,m\in\mathbb N$), assuming RH and the asymptotic formula for twisted and shifted integral moments of zeta. Under the same hypotheses, we also study a shifted case, looking at the measure $|\zeta(1/2+it+i\alpha)|^{2k}dt$, with $\alpha\in(-1,1)$. Finally we prove unconditionally the analogue result in the random matrix theory context.\end{abstract}

\section{Introduction and statement of the main results}

In 1946 Selberg \cite{Sel} proved a central limit theorem for the real part of the logarithm of the Riemann zeta function on the critical line, showing that the distribution of $\log|\zeta(1/2+it)|$ is approximately Gaussian with mean $0$ and variance $\frac{1}{2}\log\log T$, i.e. \begin{equation}\label{CLT} \frac{1}{T} \meas \Bigg \{t\in[T,2T]:\frac{ \log |\zeta (\frac{1}{2}+it)|}{\sqrt{\frac{1}{2}\log\log T}} \geq V  \Bigg \}\sim \int_V^{\infty} e^{-\frac{x^2}{2}}\frac{dx}{\sqrt{2\pi}} \end{equation}
for any fixed $V\in\mathbb R$, as $T$ goes to infinity. 
The author \cite{1.} studied the distribution of $\log|\zeta(1/2+it)|$ with respect to the weighted measure $|\zeta(1/2+it)|^2dt$ and, assuming the Riemann Hypothesis (RH), proved that it is asymptotically Gaussian with mean $\log\log T$ and variance $\frac{1}{2}\log\log T$. 
In this paper we investigate the value distribution of $\log|\zeta(1/2+it)|$ with respect to the measure 
\begin{equation}\label{wmk}|\zeta^{(m)}(1/2+it)|^{2k}dt \end{equation}
for any fixed $m,k$ non negative integers. 
The motivation is due to the study of the large values of the Riemann zeta function, i.e. the uniformity in $V$ in the central limit theorem (\ref{CLT}). In \cite{Sound}, Soundararajan speculates that an upper bound like 
$$  \frac{1}{T}\meas \Bigg \{t\in[T,2T]: \frac{\log|\zeta (1/2+it)|}{\sqrt{\frac{1}{2}\log\log T}} \geq V  \Bigg \} \ll \frac{1}{V}\exp\Big(-\frac{V^2}{2}\;\Big) $$
holds in a large range for $V$, in particular we expect (see \cite[Conjecture 2]{Radz}, see also \cite{Har1}) that for any fixed $k$
\begin{equation}\label{lvconj}  \frac{1}{T}\meas \Big \{t\in[T,2T]: \log|\zeta (1/2+it)| \geq k\log\log T  \Big \} \ll_k \frac{1}{\sqrt{\log\log T}}\frac{1}{(\log T)^{k^2}}.\end{equation}
Expressing the characteristic function of positive reals as a Mellin transform, the left hand side of (\ref{lvconj}) can be written as
$$ \frac{1}{(\log T)^{k^2}}\frac{1}{2\pi}\int_{-\infty}^{+\infty}\frac{1}{T(\log T)^{k^2}}\int_{T}^{2T}e^{iu(\log |\zeta(1/2+it)|-k\log\log T)}|\zeta(1/2+it)|^{2k}dt\frac{du}{2k+iu} $$
and is then related to the distribution of $\log|\zeta(1/2+it)|$ with respect to the weighted measure (\ref{wmk}), with $m=0$.\\

In the case $k=1,2$, we can prove a central limit theorem for $\log|\zeta(1/2+it)|$ with respect to the measure $|\zeta^{(m)}(1/2+it)|^{2k}dt$, assuming RH only.

\begin{cor}\label{derivatives}
Assume the Riemann Hypothesis, let $m$ be a non negative integer and $k=1$ or $k=2$. As $t$ varies in $T\leq t\leq 2T$, the distribution of $\log|\zeta(1/2+it)|$ is asymptotically Gaussian with mean $k\log\log T$ and variance $\frac{1}{2}\log\log T$, with respect to the weighted measure $|\zeta^{(m)}(1/2+it)|^{2k}dt$.
\end{cor}

In the case $k>2$, since not even the moments of zeta are known, we cannot expect to prove unconditionally a central limit theorem. However, 
assuming the asymptotic formula \cite[Conjecture 7.1]{HY} suggested by the recipe \cite{CFKRS} for the twisted and shifted $2k$-th moments of the Riemann zeta function 
\begin{equation}\notag 
\int_T^{2T}\Big(\frac{a}{b}\Big)^{-it}\zeta(1/2+\alpha_1+it)\cdots\zeta(1/2+\alpha_k+it)\zeta(1/2+\beta_1-it)\cdots\zeta(1/2+\beta_k-it)dt
\end{equation}
one can deal with the general case too. More precisely, we use the strategy of \cite[Theorem 1.2]{BBLR} in order to re-write Hughes and Young conjecture and we assume the following statement.

\begin{conj}[Hughes-Young]\label{HYconj}
 Let $T$ be a large parameter, $\alpha_1,\dots,\alpha_k,\beta_1,\dots,\beta_k\ll(\log T)^{-1}$, $\Phi_j$ the set of subset of $\{\alpha_1,\dots,\alpha_k\}$ of cardinality $j$, for $j=0,\dots,k$, and similarly $\Psi_j$ the set of subset of $\{\beta_1,\dots,\beta_k\}$ of cardinality $j$. If $\mathcal S \in\Phi_j$ and $\mathcal T \in \Psi_j$ then write $\mathcal S=\{\alpha_{i_1},\dots, \alpha_{i_j}\}$ and $\mathcal T=\{\beta_{l_1},\dots, \beta_{l_j}\}$ where $i_1 < i_2 < \dots < i_j$ and $l_1 < l_2 <\dots < l_j$. Let $(\alpha_{\mathcal S}; \beta_{\mathcal T})$ be the tuple obtained from $(\alpha_1,\dots ,\alpha_k;\beta_1,\dots , \beta_k)$ by replacing $\alpha_{i_r}$ with $-\beta_{l_r}$ and replacing $\beta_{l_r}$ with $-\alpha_{i_r}$
for $1 \leq r \leq j$. Consider $$ A(s)=\sum_{a\leq T^\theta}\frac{f(a)}{a^s} \hspace{0.5cm} \text{and} \hspace{0.5cm} B(s)=\sum_{b\leq T^\theta}\frac{g(b)}{b^s}$$ where $f(a)\ll a^\varepsilon$ and $g(b)\ll b^\varepsilon$ for any $\varepsilon>0$. We then conjecture that there exists a $\delta>0$, depending on $k$, such that if $\theta<\delta$ then
$$ \int_T^{2T}A(1/2+it)\overline{B(1/2+it)}\zeta(1/2+it+\alpha_1)\cdots\zeta(1/2+it+\alpha_k)\zeta(1/2-it+\beta_1)\cdots\zeta(1/2-it+\beta_{k})dt $$
equals
\begin{equation}\label{idealconj}\sum_{a,b\leq T^\theta}f(a)\overline{g(b)}\int_T^{2T}\sum_{0\leq j \leq k}\sum_{\substack{\mathcal S\in\Phi_j\\\mathcal T\in\Psi_j}}Z_{\alpha_{\mathcal S};\beta_{\mathcal T},a,b}(t)\Big(\frac{t}{2\pi}\Big)^{-\mathcal S-\mathcal T}dt+O(T^{1-\eta})\end{equation}
for some $\eta>0$, where we have written $(t/2\pi)^{-\mathcal S-\mathcal T}$ for $(t/2\pi)^{-\sum_{x\in \mathcal S}x-\sum_{y\in \mathcal T}y}$ and 
$$ Z_{\alpha;\beta,a,b}(t):=\sum_{am_1\cdots m_k=bn_1\cdots n_k}\frac{1}{\sqrt{ab}\;m_1^{1/2+\alpha_1}\cdots m_k^{1/2+\alpha_k}n_1^{1/2+\beta_1}\cdots n_k^{1/2+\beta_k}}V\bigg(\frac{m_1\cdots m_kn_1\cdots n_k}{t^k}\bigg)$$
with
$$V(x):=\frac{1}{2\pi i}\int_{(1)}\frac{G(s)}{s}(2\pi)^{-ks}x^{-s}ds$$
where $G(s)$ is an even entire function of rapid decay in any fixed strip $|\Re(s)|\leq C$ satisfying $G(0) = 1$. 
\end{conj}

With this assumption we can prove the main theorem.

\begin{thm}\label{caso2k}
Let $k,m\in\mathbb N$ and assume the Riemann Hypothesis and Conjecture \ref{HYconj} for $k$. As $t$ varies in $T\leq t\leq 2T$, the distribution of $\log|\zeta(1/2+it)|$ is asymptotically Gaussian with mean $k\log\log T$ and variance $\frac{1}{2}\log\log T$, with respect to the weighted measure $|\zeta^{(m)}(1/2+it)|^{2k}dt$.
\end{thm}

In particular, being Conjecture \ref{HYconj} known in the cases $k=1$ (see \cite{BCHB,BCR} and e.g. \cite{Bettin} for the easy modifications needed to account for the shifts) and $k=2$ (see \cite{HY} and \cite{BBLR}), we notice that Corollary \ref{derivatives} trivially follows from Theorem \ref{caso2k}. \\

We notice that Theorem \ref{caso2k} shows that the $m$-th derivative has no effect in the weighted distribution of $\log|\zeta(1/2+it)|$. This is consistent with the conjecture (see \cite[Conjecture 6.1]{Hughes} and \cite{CRS})
$$ \int_T^{2T} |\zeta(1/2+it)|^{2k-2h}|\zeta'(1/2+it)|^{2h}dt \sim c(h,k)T(\log T)^{k^2+2h}$$
which indicates that $|\zeta'(1/2+it)|^{2h}$ amplifies the contribution coming from the large values in the same way as $|\zeta(1/2+it)|^{2h}$ would do (up to a normalization of $(\log T)^{2h}$). More generally, as far as moments are concerned, the $m$-th derivative of zeta should behave like zeta itself (see e.g. \cite{CG}), up to a normalization of $\log^m T$, in accordance with Theorem \ref{caso2k}.
We also note that, while in Selberg's classical case the mean is 0 because the contribution of the small values and that of the large values of zeta balance out,  tilting with $|\zeta(1/2+it)|^{2k}$ the mean of $\log|\zeta(1/2+it)|$ moves to the right as $k$ grows and this reflects the fact that the measure $|\zeta(1/2+it)|^{2k}dt$ gives more and more weight to the large values of the Riemann zeta function.  \\

Moreover, we look at the shifted weighted measure $|\zeta(1/2+it+i\alpha)|^{2k}dt$ with $\alpha$ a real number such that $|\alpha|< 1$. As we will see, the distribution of $\log|\zeta(1/2+it)|$ is quite sensitive to the parameter $\alpha$. Indeed in computing the integral 
\begin{equation}\label{auxaux} 
\int_T^{2T}\log|\zeta(1/2+it)| |\zeta(1/2+it+i\alpha)|^{2k}dt
\end{equation}
one expects the same magnitude as in the unshifted case if $|\alpha|$ is smaller than $(\log T)^{-1}$, which is the typical scale for the Riemann zeta function. On the other hand, if $|\alpha|$ is larger then the two factors in the integral (\ref{auxaux}) start decorrelating, thus the size of the integral decreases. This phenomenon is shown in the following result, in which we use the notation $$ \widetilde\mu_\alpha=\begin{cases}\log\log T+O(1) \quad\; \text{if } |\alpha|\log T\leq 1\\-\log|\alpha|+O(1) \quad\; \text{if } |\alpha|\log T>1 \end{cases} $$ for any $\alpha\in(-1,1)$.

\begin{thm}\label{thmshift}
Let $k\in\mathbb N$ and assume the Riemann Hypothesis and Conjecture \ref{HYconj} for $k$. As $t$ varies in $T\leq t \leq 2T$, for any fixed and real $\alpha$ such that $|\alpha|< 1$, the distribution of $\log|\zeta(1/2+it)|$ is asymptotically Gaussian with mean $k\widetilde\mu_\alpha$ and variance $\frac{1}{2}\log\log T$, with respect to the measure $|\zeta(1/2+it+i\alpha)|^{2k}dt$.
\end{thm}

This theorem shows that the shift has no effect if it is smaller than $(\log T)^{-1}$. On the contrary, for larger values of the shift the mean gets smaller, for instance if $$\alpha = \frac{(\log T)^\delta}{\log T}$$ with $\delta\in (0,1)$ then the mean is $\sim(1-\delta)k\log\log T$. In this shifted case too, if $k\leq 2$ then Theorem \ref{thmshift} holds assuming RH only.\\

Lastly, we show that in the random matrix theory setting, an analogous weighted central limit theorem can be proved unconditionally. We consider the characteristic polynomials
$$ Z=Z(U,\theta)=\det(I-Ue^{-i\theta})$$
of $N\times N$ unitary matrices $U$ and we investigate their distribution of values with respect to the circular unitary ensemble (CUE). It has been conjectured that the limiting distribution of the non-trivial zeros of the Riemann zeta function, on the scale of their mean spacing, is the same as that of the eigenphases $\theta_n$ of matrices in the CUE in the limit as $N\to\infty$ (see e.g. \cite{RM2000}). Then we consider a tilted version of the Haar measure and we have the following theorem.
\begin{thm}\label{rmt}
As $N\to \infty$, the value distribution of $\log |Z|$ is asymptotically Gaussian with mean $k\log N$ and variance $\frac{1}{2}\log N$ with respect to the measure $|Z|^{2k}d_{Haar}$ .
\end{thm}
As usual, the correspondence with Theorem \ref{caso2k} holds if we identify the mean density of the eigenangles $\theta_n$, $N/2\pi$, with the mean density of the Riemann zeros at a height $T$ up the critical line, $\frac{1}{2\pi}\log\frac{T}{2\pi}$, i.e. if $$ N=\log\frac{T}{2\pi}. $$


\section{Proof of Theorems \ref{caso2k} and \ref{thmshift}} 

To prove both the theorems, we introduce a set of shifts $\alpha_1,\dots,\alpha_k,\beta_1,\dots,\beta_k$ and we denote for the sake of brevity
$$ \zeta_{\alpha,\beta}(t) := \zeta(1/2+\alpha_1+it)\cdots \zeta(1/2+\alpha_k+it)\zeta(1/2+\beta_1-it)\cdots\zeta(1/2+\beta_k-it).$$ 
The general strategy of the proof is similar to the one of Theorem 1 in \cite{1.}, but here we avoid the detailed combinatorial analysis we performed in \cite{1.}, by working with Euler products instead of Dirichlet series, inspired by \cite[Proposition 5.1]{ABBRS}. The first step is then approximating the logarithm of the Riemann zeta function with a suitable Dirichlet polynomial. Let's denote 
\begin{equation}\label{poly1}\widetilde P(t) :=\sum_{p\leq x}\frac{1}{p^{1/2+it}}\end{equation} 
where $x :=T^\varepsilon$, with $\varepsilon := (\log \log \log T)^{-1}$. Now we show that $\log|\zeta(1/2+it)|$ has the same distribution as $\Re \widetilde P(t)$ with respect to the measure $\zeta_{\alpha,\beta}(t)dt$. This is achieved in the following proposition by bounding the second moment of the difference. 

\begin{prop}\label{prop1}
For $k$ a non negative integer, assume Conjecture \ref{HYconj} and RH. Let T be a large parameter and $\alpha_1,\dots,\alpha_k$, $\beta_1,\dots,\beta_k\in\mathbb C$ such that $|\alpha_i|,|\beta_j|< 1$, $|\Re(\alpha_i)|,|\Re(\beta_j)|\leq (\log T)^{-1}$ and $|\alpha_i-\beta_j|\ll(\log T)^{-1}$ for all $1\leq i,j\leq k$. Then we have:
$$\int_T^{2T}\big|\log|\zeta(1/2+it)|-\Re \widetilde P(t)\big|^2\zeta_{\alpha,\beta}(t)dt\ll_k T(\log T)^{k^2}(\log\log\log T)^{5/2}.$$
\end{prop}

\proof
From Tsang's work \cite[equation (5.15)]{Tsang} we know that 
\begin{equation}\label{approx1} \log \zeta(1/2+it)-\widetilde P(t) = S_1 +S_2 +S_3 + O(R) - L(t)  \end{equation}
with 
$$ S_1:= \sum_{p\leq x}\big ( p^{-1/2-4/\log x}-p^{-1/2}\big )p^{-it},\hspace{1cm} 
S_2:=\sum_{\substack{p^r\leq x \\ r\geq 2}}\frac{p^{-r(1/2+4/\log x+it)}}{r}$$
$$ S_3:=\sum_{x<n\leq x^3}\frac{\Lambda(n)}{\log n}n^{-1/2-4/\log x-it}, \hspace{1cm}
R:= \frac{5}{\log x}\bigg ( \bigg | \sum_{n\leq x^3}\frac{\Lambda (n)}{n^{1/2+4/\log x+it}} \bigg |+\log T \bigg ) $$
$$ L(t):=\sum_\rho \int_{1/2}^{1/2+4/\log x}\left ( \frac{1}{2}+\frac{4}{\log x}-u \right )\frac{1}{u+it-\rho}\frac{1}{\frac{1}{2}+\frac{4}{\log x}-\rho}du,$$
where the sum in the definition of $L(t)$ is over all the non-trivial zeros of $\zeta$, then we have to bound the second moments of the terms on the right hand side with respect to weighted measure $\zeta_{\alpha,\beta}(t)dt$, by using Conjecture \ref{HYconj} (note that we are allowed to apply the conjecture, since the shifts are small up to a change of variable, being $|\alpha_i-\beta_j|\ll(\log T)^{-1}$).\\ 

Let's start with $\int_T^{2T}|S_1|^2\zeta_{\alpha,\beta}(t)dt$ which is 
\begin{equation}\begin{split}\label{S1.1}
&\ll\bigg|\int_T^{2T}\sum_{0\leq j\leq k}\sum_{\substack{\mathcal S\in\Phi_j \\ \mathcal T\in \Psi_j}}\Big(\frac{t}{2\pi}\Big)^{-\mathcal S-\mathcal T}\sum_{p,q\leq x}(p^{-4/\log x}-1)(q^{-4/\log x}-1)Z_{\alpha_{\mathcal S},\beta_{\mathcal T},p,q}(t)dt\bigg|\\&
\ll \int_T^{2T}\sum_{0\leq j\leq k}\sum_{\substack{\mathcal S\in\Phi_j \\ \mathcal T\in \Psi_j}}\bigg|\sum_{a,b\leq x}\mathbb 1_{x}(a)\mathbb 1_{x}(b)(a^{-\frac{4}{\log x}}-1)(b^{-\frac{4}{\log x}}-1)Z_{\alpha_{\mathcal S},\beta_{\mathcal T},a,b}(t)\bigg|dt
\end{split}\end{equation}
where $\mathbb 1_x(\cdot)$ is the indicator function of primes up to $x$. If we denote $\Omega(n)$ the function which counts the number of prime factors of $n$ with multiplicity, then we have 
\begin{equation}
\label{moltiplicativizzare}\sum_{n\leq x}\mathbb 1_{x}(n)f(n)=
\sum_{\substack{n:\\p|n\implies p\leq x}}\partial_z\big[z^{\Omega(n)}f(n)\big]_{z=0}\end{equation}
hence we get
\begin{equation}\label{S1.2}
\int_T^{2T}|S_1|^2\zeta_{\alpha,\beta}(t)dt
\ll \int_T^{2T}\sum_{0\leq j\leq k}\sum_{\substack{\mathcal S\in\Phi_j \\ \mathcal T\in \Psi_j}}\bigg|\partial_z\partial_w \bigg[\frac{1}{2\pi i}\int_{(1)}\frac{G(s)}{s}\Big(\frac{t}{2\pi}\Big)^{ks}I_{\alpha,\beta}^{j,\mathcal S,\mathcal T}(z,w;s)ds\bigg]_{\substack{z=0\\w=0}}\bigg|dt
\end{equation}
where
\begin{equation}\label{S1.3}
I_{\alpha,\beta}^{j,\mathcal S,\mathcal T}(z,w;s):=\sum_{\substack{a,b: \\ p|ab\implies p\leq x}} z^{\Omega(a)}w^{\Omega(b)}g_x(a)g_x(b)\widetilde Z_{\alpha_{\mathcal S},\beta_{\mathcal T},a,b}(s)
\end{equation}
with $g_x$ the multiplicative function defined by $g_x(p^\alpha)=p^{-4\alpha/\log x}-1$ and
\begin{equation}\label{widetildeZ}\widetilde Z_{\alpha,\beta,a,b}(s):=\sum_{am_1\cdots m_k=bn_1\cdots n_k}\frac{1}{\sqrt{ab}\;m_1^{1/2+\alpha_1+s}\cdots m_k^{1/2+\alpha_k+s}n_1^{1/2+\beta_1+s}\cdots n_k^{1/2+\beta_k+s}} .\end{equation}
We now analyze the first term $I_{\alpha,\beta}^{0}(z,w;s)$ and then we will see how to apply the method to deal with all the others. 
Since the sum in the definition (\ref{S1.3}) is multiplicative we have
\begin{equation}\begin{split}\label{S1.4}
I_{\alpha,\beta}^{0}(z,w;s)&=
\sum_{\substack{am_1\cdots m_k=bn_1\cdots n_k: \\ p|ab\implies p\leq x}}\frac{z^{\Omega(a)}w^{\Omega(b)}g_x(a)g_x(b)}{\sqrt{ab}\; m_1^{\frac{1}{2}+\alpha_1+s}\cdots m_k^{\frac{1}{2}+\alpha_k+s}n_1^{\frac{1}{2}+\beta_1+s}\cdots n_k^{\frac{1}{2}+\beta_k+s}}\\&
=\prod_{p\leq x}\;\sum_{\substack{a+m_1+\dots+m_k= \\ =b+n_1+\dots+n_k}}\frac{z^{\Omega(p^a)}w^{\Omega(p^b)}g_x(p^a)g_x(p^b)}{p^{\frac{a}{2}+\frac{b}{2}+m_1(\frac{1}{2}+\alpha_1+s)+\dots+m_k(\frac{1}{2}+\alpha_k+s)+n_1(\frac{1}{2}+\beta_1+s)+\dots+ n_k(\frac{1}{2}+\beta_k+s)}}\\&
\hspace{1.4cm}\cdot\prod_{p> x}\;\sum_{\substack{m_1+\dots+m_k= \\=n_1+\dots+n_k}}\frac{1}{p^{m_1(\frac{1}{2}+\alpha_1+s)+\dots+m_k(\frac{1}{2}+\alpha_k+s)+m_1(\frac{1}{2}+\beta_1+s)+\dots + n_k(\frac{1}{2}+\beta_k+s)}}
\end{split}\end{equation}
and by putting in evidence the first terms in the Euler products, this is
\begin{equation}\begin{split}\label{S1.5}
&=A_{\alpha,\beta}(z,w;s) \prod_{p\leq x}\bigg(1+\frac{zwg_x(p)^2}{p}+\frac{zg_x(p)}{p^{1+\beta_1+s}}+\dots+\frac{wg_x(p)}{p^{1+\alpha_k+s}}\bigg)\\&
\hspace{0.5cm}\cdot\prod_{p\leq x}\bigg(1+\frac{1}{p^{1+\alpha_1+\beta_1+2s}}+\dots+\frac{1}{p^{1+\alpha_k+\beta_k+2s}}\bigg)\prod_{p>x}\bigg(1+\frac{1}{p^{1+\alpha_1+\beta_1+2s}}+\dots+\frac{1}{p^{1+\alpha_k+\beta_k+2s}}\bigg)\\&
=A_{\alpha,\beta}^*(z,w;s) \prod_{i,j=1}^k\zeta(1+\alpha_i+\beta_j+2s)\cdot \exp\bigg(\sum_{p\leq x}\bigg\{\frac{zwg_x(p)^2}{p}+\frac{zg_x(p)}{p^{1+\beta_1+s}}+\dots+\frac{wg_x(p)}{p^{1+\alpha_k+s}}\bigg\}\bigg) 
\end{split}\end{equation}
where $A_{\alpha,\beta}(z,w;s)$ and $A_{\alpha,\beta}^*(z,w;s) $  are arithmetical factors (Euler products) converging absolutely in a half-plane $\Re(s)>-\delta$ for some $\delta>0$, uniformly for $|z|,|w|\leq 1$, such that their derivatives with respect to $z$ and $w$ at 0 also converge in the same half plane.
We now have extracted the polar part, hence we are ready to shift the integral over $s$ in (\ref{S1.2}) to the left of zero. To do so, it is convenient to prescribe the same conditions as in \cite[Remarks after Lemma 2.1]{BBLR}, assuming that $G(s)$ vanishes at $s=-\frac{\alpha_i+\beta_j}{2}$ for all $i,j$, so that the only pole we pick in the contour shift is at $s=0$. Moreover we assume that the shifts are such that $|\alpha_i+\beta_j|\gg(\log T)^{-1}$ for every $1\leq i,j\leq k$, so that
\begin{equation}\label{S1.4,5} \prod_{i,j=1}^k|\zeta(1+\alpha_i+\beta_j)|\ll_k(\log T)^{k^2}. \end{equation}
Hence by (\ref{S1.4}), (\ref{S1.5}) and (\ref{S1.4,5}), we get
\begin{equation}\begin{split}\notag
\partial_z\partial_w\bigg[\frac{1}{2\pi i}&\int_{(1)}\frac{G(s)}{s}\Big(\frac{t}{2\pi}\Big)^{ks}I_{\alpha,\beta}^{0}(z,w;s)ds\bigg]_{\substack{z=0\\w=0}} 
\\& \ll_k (\log T)^{k^2} \bigg|\partial_z\partial_w\bigg[ \exp\bigg(\sum_{p\leq x}\bigg\{\frac{zwg_x(p)^2}{p}+\frac{zg_x(p)}{p^{1+\beta_1}}+\dots+\frac{wg_x(p)}{p^{1+\alpha_k}}\bigg\}\bigg)\bigg]_{\substack{z=0\\w=0}} \bigg|
\\& \ll_k (\log T)^{k^2}\bigg(\sum_{p\leq x}\frac{|g_x(p)|}{p}\bigg)^2 \ll_k (\log T)^{k^2} (\log\log\log T)^2
\end{split}\end{equation}
being $|\Re(\alpha_i)|,|\Re(\beta_j)|\leq (\log T)^{-1}$ and $\sum_{p\leq x}|p^{-4/\log x}-1|/p \ll 1$; this last bound can be obtained by using the Taylor's approximation $e^{-z}=1+O(z)$ for $z\ll1$, which yields 
$$ \sum_{p\leq x}\frac{|p^{-4/\log x}-1|}{p}= \sum_{p\leq x}\frac{|e^{-4\log p/\log x}-1|}{p}\ll \frac{1}{\log x}\sum_{p\leq x}\frac{\log p}{p}\ll 1$$ by Mertens' first theorem.
All the other terms $I_{\alpha,\beta}^{j,\mathcal S,\mathcal T}(z,w)$ can be treated exactly in the same way as $I_{\alpha,\beta}^0(z,w)$ by assuming that $|\alpha_i\pm\beta_j|\gg (\log T)^{-1}$ for every $i,j$, since they only differ from the first case by permutations and changes of signs of the shifts. Therefore we get 
\begin{equation}\begin{split}\notag
\partial_z\partial_w\bigg[\frac{1}{2\pi i}&\int_{(1)}\frac{G(s)}{s}\Big(\frac{t}{2\pi}\Big)^{ks}I_{\alpha,\beta}^{j,\mathcal S,\mathcal T}(z,w;s)ds\bigg]_{\substack{z=0\\w=0}}  
\ll_k(\log T)^{k^2} (\log\log\log T)^2
\end{split}\end{equation}
provided that 
\begin{equation}\begin{split}\label{S1.7}
|\alpha_i\pm\beta_j|\gg (\log T)^{-1}\quad \text{for every } 1\leq i,j\leq k.
\end{split}\end{equation}
Plugging this into (\ref{S1.2}), we prove that 
\begin{equation}\begin{split}\label{S1.8}
\int_T^{2T}|S_1|^2\zeta_{\alpha,\beta}(t)dt\ll_k T(\log T)^{k^2}(\log\log\log T)^2.
\end{split}\end{equation}
Moreover, since the left hand side in (\ref{S1.8}) is holomorphic in terms of the shifts, although we have proved the above for $\alpha_i,\beta_j$ such that (\ref{S1.7}) holds, the maximum modulus principle can be applied to obtain the bound to the enlarged domain we need.

We will treat similarly also the other pieces. As regards the second one, we have that $\int_T^{2T}|S_2|^2\zeta_{\alpha,\beta}(t)dt$ is
\begin{equation}\begin{split}\notag
\ll  \int_T^{2T}\sum_{0\leq j\leq k}\sum_{\substack{\mathcal S\in\Phi_j \\ \mathcal T\in \Psi_j}}\bigg| \sum_{\substack{p_1^{r_1},p_2^{r_2}\leq x \\ r_1,r_2\geq 2}}\frac{p_1^{-\frac{4r_1}{\log x}}}{r_1}\frac{p_2^{-\frac{4r_2}{\log x}}}{r_2}Z_{\alpha_{\mathcal S},\beta_{\mathcal T},p_1^{r_1},p_2^{r_2}}(t)\bigg|dt.
\end{split}\end{equation}
As before, we analyze the first term only since all the others are completely analogous. The term for $j=0$ is 
\begin{equation}\begin{split}\notag
\frac{1}{2\pi i}\int_{(1)}\frac{G(s)}{s}\Big(\frac{t}{2\pi}\Big)^{ks}\sum_{\substack{p_1^{r_1}m_1\cdots m_k =p_2^{r_2}n_1\cdots n_k : \\ p_1^{r_1},p_2^{r_2}\leq x, \; r_1,r_2\geq 2}}&\frac{p_1^{-\frac{4r_1}{\log x}}p_2^{-\frac{4r_2}{\log x}}}{r_1r_2}\frac{1}{\sqrt{p_1^{r_1}p_2^{r_2}}\; m_1^{\frac{1}{2}+\alpha_1+s}\cdots n_k^{\frac{1}{2}+\beta_k+s}}ds
\end{split}\end{equation}
and this time, because of the condition $r_1,r_2\geq 2$, when we estimate the sum via the first terms of its Euler product we just get that the above is
$$ \ll_k \prod_{i,j=1}^k|\zeta(1+\alpha_i+\beta_j)| $$
and, applying the same machinery as before, this yields
\begin{equation}\label{S2.1} 
\int_T^{2T}|S_2|^2\zeta_{\alpha,\beta}(t)dt \ll_k T(\log T)^{k^2}. 
\end{equation}

We use the same approach in order to bound the second moment of $S_3$ as well, which is 
\begin{equation}\begin{split}\notag
&\ll_k \int_T^{2T} \bigg| \frac{1}{2\pi i}\int_{(1)}\frac{G(s)}{s}\Big(\frac{t}{2\pi}\Big)^{ks}\sum_{\substack{am_1\cdots m_k=bn_1\cdots n_k \\ x<a,b<x^3}}\frac{\Lambda(a)\Lambda(b)}{\log a\log b}\frac{a^{-4/\log x}\;b^{-4/\log x}}{\sqrt{ab}\;m_1^{\frac{1}{2}+\alpha_1+s}\cdots n_k^{\frac{1}{2}+\beta_k+s}} \bigg|dt\\&
\ll_k T\bigg(\sum_{x<p\leq x^3}\frac{1}{p}\bigg)^2\prod_{i,j=1}^k|\zeta(1+\alpha_i+\beta_j)| \ll_k T(\log T)^{k^2}
\end{split}\end{equation}
since the sum $\sum_{x<p\leq x^3}p^{-1}$ is bounded.\\

We deal with $R$ in the same way and we get that
\begin{equation}\label{R.1}
\int_T^{2T}|R|^2\zeta_{\alpha,\beta}(t)dt \ll_k T(\log T)^{k^2}(\log\log\log T)^2 \end{equation}
where the extra factor with the triple $\log$ comes from the second term in the definition of $R$, while the first one can be treated analogously to $S_3$.\\

Finally we have to bound the second moment of $\Re L(t)$. To do so, in view of equations (2.8) and (2.9) in \cite{1.}, it suffices to study
\begin{equation}\begin{split}\label{L.1}
\frac{1}{(\log x)^{2}}\int_T^{2T} \Big (\log ^+ \frac{1}{\eta _t \log x} \Big )^{2}\Big ( \Big | \sum_{n\leq x^3} \frac{\Lambda(n)n^{-4/\log x}}{n^{1/2+it}}\Big |+\log T\Big )^{2}\zeta_{\alpha,\beta}(t)dt  
\end{split}\end{equation}
where $\eta_t:=\min_\rho|t-\gamma|$ and $\log ^+ t:=\max(\log t,0)$, with the aim of proving that this is $\ll_k T(\log T)^{k^2}(\log\log\log T)^{5/2}$. By applying the Cauchy-Schwarz inequality, the above is
\begin{equation}\begin{split}\label{L.2}
\leq& \bigg(\int_T^{2T} \Big (\frac{1}{\log x} \Big | \sum_{n\leq x^3} \frac{\Lambda(n)n^{-4/\log x}}{n^{1/2+it}}\Big |+\frac{1}{\varepsilon}\Big )^{4}|\zeta(1/2+\beta_1-it)|^2\cdots |\zeta(1/2+\beta_k-it)|^2dt \bigg)^{\frac{1}{2}}\\&
\hspace{3.1cm} \cdot\bigg( \int_T^{2T} \Big (\log ^+ \frac{1}{\eta _t \log x} \Big )^{4}|\zeta(1/2+\alpha_1+it)|^2\cdots |\zeta(1/2+\alpha_k+it)|^2dt  \bigg)^{\frac{1}{2}}
\end{split}\end{equation}
and the first term can be treated as $R$ above and it is $\ll_k \sqrt{T(\log T)^{k^2}(\log\log\log T)^4}$. 
We now conclude the proof, bounding the second term in (\ref{L.2}). If we denote $\tau:=[T-\frac{1}{\log x},2T+\frac{1}{\log x}]$ then we have
\begin{equation}\begin{split}\label{L.3}
\int_T^{2T} \Big (&\log ^+ \frac{1}{\eta _t \log x}\Big )^{4}  |\zeta(1/2+\alpha_1+it)|^2\cdots |\zeta(1/2+\alpha_k+it)|^2dt \\ &
\leq \sum_{\gamma\in\tau} \int_{-1/\log x}^{1/\log x} \Big (\log ^+ \frac{1}{|w| \log x}\Big )^{4} \prod_{j=1}^k\big |\zeta\big (1/2+\alpha_j+i(w+\gamma)\big )\big |^2dw \\ & 
= \sum_{\gamma\in\tau} \int_{-1}^1 \Big (\log ^+ \frac{1}{|t|}\Big )^{4}\prod_{j=1}^k\bigg |\zeta\Big (1/2+\alpha_j+i\Big(\gamma+\frac{t}{\log x}\Big)\Big )\bigg |^2\frac{dt}{\log x}\\&
\ll\frac{1}{\log x} \int_{-1}^1 (\log |t|)^{4}\prod_{j=1}^k\bigg(\sum_{\gamma\in\tau} \Big |\zeta\Big (1/2+i\gamma+\alpha_j+i\frac{t}{\log x}\Big |^{2k}\bigg)^{1/k}dt
\end{split}\end{equation}
by H\"older inequality. The remaining sum can be bounded under RH in view of Kirila's \cite[Theorem 1.2]{Kirila} and Milinovich's \cite{Milinovich} works, which generalize the well known result due to Gonek about sum over the non-trivial zeros of zeta \cite[Corollary 1]{Gonek}. Indeed, since the shifts $\alpha_j+i\frac{t}{\log x}$ in the sum over zeros in (\ref{L.3}) have modulus $\leq 1$ and real part $\leq (\log T)^{-1}$ in absolute value, then we have 
$$ \sum_{\gamma\in\tau}\Big |\zeta\Big (1/2+i\gamma+\alpha_j+i\frac{t}{\log x}\Big )\Big |^{2k}\ll_k T\log T(\log T)^{k^2} $$ 
 for every $j=1,\dots,k$ and putting this into (\ref{L.3}) we get
\begin{equation}\label{L.4}
\int_T^{2T} \Big (\log ^+ \frac{1}{\eta _t \log x}\Big )^{4}  |\zeta(1/2+\alpha_1+it)|^2\cdots |\zeta(1/2+\alpha_k+it)|^2dt \ll_k \frac{T\log T}{\log x}(\log T)^{k^2}.
\end{equation}
Plugging (\ref{L.4}) into (\ref{L.2}) we prove that
\begin{equation}\begin{split}\label{L.5}
\int_T^{2T} |\Re L(t)|^2\zeta_{\alpha,\beta}(t)dt &
\ll_k\sqrt{T(\log T)^{k^2}(\log\log\log T)^4}\sqrt{T(\log T)^{k^2}\log\log\log T}\\&
\ll_k T(\log T)^{k^2}(\log\log\log T)^{5/2}
\end{split}\end{equation}
concluding the proof of the proposition.
\endproof

The second step is getting rid of the small primes, showing that their contribution does not affect the distribution asymptotically. This simple fact will simplify the third and last step of the proof, as we will see in the following. Let's define 
\begin{equation}\label{poly}P(t) :=\sum_{p\in X}\frac{1}{p^{1/2+it}}\end{equation}
where $X:=(\log T,x]$ (we recall that $x=T^\varepsilon$ and $\varepsilon=(\log\log\log T)^{-1}$).

\begin{prop}\label{prop2}
For $k$ a non negative integer, assume Conjecture \ref{HYconj} and RH. Let T be a large parameter and $\alpha_1,\dots,\alpha_k$, $\beta_1,\dots,\beta_k\in\mathbb C$ such that $|\alpha_i|,|\beta_j|< 1$, $|\Re(\alpha_i)|,|\Re(\beta_j)|\leq (\log T)^{-1}$ and $|\alpha_i-\beta_j|\ll(\log T)^{-1}$ for all $1\leq i,j\leq k$. Then we have:
$$\int_T^{2T}\big|\Re\widetilde P(t)-\Re  P(t)\big|^2\zeta_{\alpha,\beta}(t)dt\ll_k T(\log T)^{k^2}(\log\log\log T)^2.$$ 
\end{prop}

\proof
This can be proved with the same method used in Proposition \ref{prop1}. We recall that $\mathbb 1_{\log T}(\cdot)$ denotes the indicator function of primes up to $\log T$. We start by studying
\begin{equation}\label{P2.10}
\frac{1}{2\pi i}\int_{(1)}\frac{G(s)}{s}\Big(\frac{t}{2\pi}\Big)^{ks}\sum_{\substack{am_1\cdots m_k= \\ =bn_1\cdots n_k}}\frac{\mathbb 1_{\log T}(a)\mathbb 1_{\log T}(b)}{\sqrt{ab}\; m_1^{\frac{1}{2}+\alpha_1+s}\cdots n_k^{\frac{1}{2}+\beta_k+s}}ds
\end{equation}
and, as usual, we estimate the sum with the first terms of its Euler product, extract the polar part, shift the integral over $s$ to the left, getting that (\ref{P2.10}) is
\begin{equation}\notag
\ll_k \prod_{i,j=1}^k|\zeta(1+\alpha_i+\beta_j)|\cdot\bigg|\partial_z\partial_w\bigg[
\exp\bigg(\sum_{p\leq \log T}\bigg\{\frac{zw}{p}+\frac{z}{p^{1+\beta_1}}+\dots+\frac{w}{p^{1+\alpha_k}}\bigg\}\bigg)
\bigg]_{\substack{z=0\\w=0}}\bigg|
\end{equation}
and by the same argument as in the proof of Proposition \ref{prop1} the above is
\begin{equation}\notag
\ll_k (\log T)^{k^2}\bigg(\sum_{p\leq \log T}\frac{1}{p}\bigg)^2\ll (\log T)^{k^2} (\log\log\log T)^2
\end{equation}
and this concludes the proof.
\endproof

Finally we investigate the distribution of the polynomial $\Re P(t)$, which has the same distribution as $\log|\zeta(1/2+it)|$ thanks to Propositions \ref{prop1} and \ref{prop2}. The most natural method to do so is studying the moments and this is achieved in the following result.

\begin{prop}\label{prop3}
For $k$ a non negative integer, assume Conjecture \ref{HYconj} and RH. Let T be a large parameter and $\alpha_1,\dots,\alpha_k$, $\beta_1,\dots,\beta_k\in\mathbb C$ such that $|\alpha_i|,|\beta_j|< 1$, $|\Re(\alpha_i)|,|\Re(\beta_j)|\leq (\log T)^{-1}$ and $|\alpha_i-\beta_j|\ll(\log T)^{-1}$ for all $1\leq i,j\leq k$. Denote $\mathcal L:=\sum_{p\in X}\frac{1}{p}\sim\log\log T$ and $\mu\in \mathbb R$ such that $\mu\ll\log\log T$. Then for every fixed integer $n$ we have
\begin{equation}\begin{split}\notag
&\int_T^{2T}(\Re P(t)-k\mu)^n\zeta_{\alpha,\beta}(t)dt \\&\hspace{1.4cm}
= \int_T^{2T}\sum_{0\leq j\leq k}\sum_{\substack{\mathcal S\in\Phi_j \\ \mathcal T\in\Psi_j}}\Big(\frac{t}{2\pi}\Big)^{-\mathcal S-\mathcal T}M_{\alpha_{\mathcal S},\beta_{\mathcal T}}(0)\;
\partial_z^n\bigg[e^{\frac{z^2}{4}\mathcal L-kz\mu}\exp\bigg(\frac{z}{2}\sum_{p\in X}\frac{g_p(\mathcal S,\mathcal T)}{p}\bigg)\bigg]_{z=0}dt \\&\hspace{1.4cm} \quad + O_{k,n}\Big(T(\log T)^{k^2-1+\varepsilon}\Big)
\end{split}\end{equation}
where 
$$g_p(\mathcal S,\mathcal T):=\sum_{x_1\not\in\mathcal S}p^{-x_1}+\sum_{x_2\not\in\mathcal T}p^{-x_2}+\sum_{x_3\in\mathcal S}p^{x_3}+\sum_{x_4\in\mathcal T}p^{x_4}$$
and
$$M_{\alpha,\beta}(s):=\sum_{m_1\cdots m_k=n_1\cdots n_k}\frac{1}{m_1^{\frac{1}{2}+\alpha_1+s}\cdots m_k^{\frac{1}{2}+\alpha_k+s}n_1^{\frac{1}{2}+\beta_1+s}\cdots n_k^{\frac{1}{2}+\beta_k+s}}$$
so that $M_{\alpha,\beta}(0)$ is the first term of the moment of $\zeta_{\alpha,\beta}$, predicted by the recipe \cite{CFKRS} (the sum which defines $M_{\alpha,\beta}(s)$ does not converge for $s=0$, so we appeal to \cite[Theorem 2.4.1]{CFKRS} for the analytic continuation).
\end{prop}

\proof
Expanding out the powers, since $2\Re(z)=z+\overline z$, we get
\begin{equation}\label{P3.1}
\int_T^{2T}(\Re P(t)-k\mu)^n\zeta_{\alpha,\beta}(t)dt=
\sum_{j+h=n}\binom{n}{h}(-k\mu)^j2^{-h}\sum_{r+s=h}\binom{h}{r}\int_T^{2T}P(t)^r\overline{P(t)}^s\zeta_{\alpha,\beta}(t)dt
\end{equation}
and using Conjecture \ref{HYconj}, ignoring the error term which is negligible in this context, the inner integral in (\ref{P3.1}) is
\begin{equation}\label{P3.2}
\int_T^{2T}\sum_{0\leq j\leq k}\sum_{\substack{\mathcal S\in\Phi_j \\ \mathcal T\in\Psi_j}}\Big(\frac{t}{2\pi}\Big)^{-\mathcal S-\mathcal T} \sum_{a,b}\mathbb 1_X^{*r}(a)\mathbb 1_X^{*s}(b)Z_{\alpha_{\mathcal S},\beta_{\mathcal T},a,b}(t) dt
\end{equation}
where $\mathbb 1_X^{*r}(a)$ denotes the indicator function of primes in the interval $X$, self-convoluted $r$ times. If we define temporarily the multiplicative function $g$ given by $g(p^n) = 1/n!$, then the inner sum over $a,b$ in (\ref{P3.2}) equals
\begin{equation}\label{P3.3}
\partial_z^r\partial_w^s\bigg[ \sum_{\substack{a,b:\\p|ab\implies p\in X}} z^{\Omega(a)}w^{\Omega(b)}g(a)g(b) Z_{\alpha_{\mathcal S},\beta_{\mathcal T},a,b}(t)\bigg]_{\substack{z=0\\w=0}}
\end{equation}
then plugging (\ref{P3.2}) and (\ref{P3.3}) into (\ref{P3.1}), recollecting together the powers we expanded before, we get
\begin{equation}\label{P3.4}
\int_T^{2T}\frac{1}{2\pi i}\int_{(1)}\frac{G(s)}{s}\Big(\frac{t}{2\pi}\Big)^{ks}\sum_{0\leq j\leq k}\sum_{\substack{\mathcal S\in\Phi_j \\ \mathcal T\in\Psi_j}}\Big(\frac{t}{2\pi}\Big)^{-\mathcal S-\mathcal T} \partial_z^n\Big[e^{-k\mu z}I_{\alpha,\beta}^{j,\mathcal S,\mathcal T}(z;s)\Big]_{z=0}dsdt.
\end{equation}
with
\begin{equation}\label{P3.5}
I_{\alpha,\beta}^{j,\mathcal S,\mathcal T}(z;s):=\sum_{\substack{a,b: \\ p|ab\implies p\in X}}\Big(\frac{z}{2}\Big)^{\Omega(a)+\Omega(b)}g(a)g(b)\widetilde Z_{\alpha_{\mathcal S},\beta_{\mathcal T},a,b}(s)
\end{equation}
(see (\ref{widetildeZ}) for the definition of $\widetilde Z_{\alpha,\beta,a,b}(s)$). Now study the first term in (\ref{P3.4}), i.e. $j=0$, since all the other terms can be understood from the first one with a slight modification. Then we look at
\begin{equation}\begin{split}\label{P3.6}
&I_{\alpha,\beta}^{0}(z;s)=\sum_{\substack{am_1\cdots m_k=bn_1\cdots n_k \\ p|ab\implies p\in X}}\frac{(z/2)^{\Omega(a)+\Omega(b)}g(a)g(b)}{\sqrt{ab}\;m_1^{\frac{1}{2}+\alpha_1+s}\cdots m_k^{\frac{1}{2}+\alpha_k+s}n_1^{\frac{1}{2}+\beta_1+s}\cdots n_k^{\frac{1}{2}+\beta_k+s}}\\&
=\prod_{p\in X}\sum_{\substack{a+m_1+\dots+m_k=\\=b+n_1+\dots+n_k}}\frac{(z/2)^{a+b}g(p^a)g(p^b)}{p^{\frac{a}{2}+\frac{b}{2}+m_1(\frac{1}{2}+\alpha_1+s)+\dots+n_k(\frac{1}{2}+\beta_k+s)}}
\prod_{p\not\in X}\sum_{\substack{m_1+\dots+m_k=\\=n_1+\dots+n_k}}\frac{1}{p^{m_1(\frac{1}{2}+\alpha_1+s)+\dots+n_k(\frac{1}{2}+\beta_k+s)}}\\&
=\exp\bigg(\sum_{p\in X}\bigg\{\frac{(z/2)^2}{p}+\frac{z/2}{p^{1+\alpha_1+s}}+\dots+\frac{z/2}{p^{1+\alpha_k+s}}+\frac{z/2}{p^{1+\beta_1+s}}+\frac{z/2}{p^{1+\beta_k+s}}\bigg\}\bigg)F_{X,\alpha,\beta}(z;s)M_{\alpha,\beta}(s)
\end{split}\end{equation}
where $F_{X,\alpha,\beta}(z;s)$ is an arithmetical factor (Euler product) converging absolutely in a product of half-planes containing the origin, such that $F_{X,\alpha,\beta}(z;0)$ is holomorphic at $z=0$, $F_{X,\alpha,\beta}(0,0)=1$ and all its derivatives at $z=0$ are small, i.e. $\partial_z^c[F_{X,\alpha,\beta}(z,0)]_{z=0}\ll (\log T)^{-1}$ for any positive integer $c$.
Now we want to shift the integral over $s$ in (\ref{P3.4}) to the left of zero, picking the contribution of the (only) pole at $s=0$. To do so, we appeal to the meromorphic continuation of the function $M_{\alpha,\beta}(s)$, see \cite[Theorem 2.4.1]{CFKRS}; thus we can shift the path of integration to the vertical line (say) $\Re(s)=-\frac{1}{10}$, where the integral is trivially bounded by $\ll T^{1-1/10+\varepsilon}$ for any positive $\varepsilon$. Moreover, the contribution from the pole at $s=0$ gives
\begin{equation}\notag
T\partial_z^n\bigg[e^{-k\mu z}e^{\frac{z^2}{4}\mathcal L}\exp\bigg(\frac{z}{2}\sum_{p\in X}\frac{p^{-\alpha_1}+\dots+p^{-\beta_k}}{p}\bigg) F_{X,\alpha,\beta}(z;0)\bigg]_{z=0}M_{\alpha,\beta}(0).
\end{equation}
Thanks to the bounds for $\mu$, $\mathcal L$ and the derivatives of $F_{X,\alpha,\beta}$, being $F_{X,\alpha,\beta}(0;0)=1$ and $M_{\alpha,\beta}(0)\ll(\log T)^{k^2}$ (again this is due to a similar argument as in the proof of Proposition \ref{prop1}, when we assume extra conditions on the shifts and then we appeal to the maximum modulus principle), we get that the term for $j=0$ in (\ref{P3.4}) equals
\begin{equation}\notag
T\partial_z^n\bigg[e^{-k\mu z}e^{\frac{z^2}{4}\mathcal L}\exp\bigg(\frac{z}{2}\sum_{p\in X}\frac{p^{-\alpha_1}+\dots+p^{-\beta_k}}{p}\bigg) \bigg]_{z=0}M_{\alpha,\beta}(0) + O_{k,n}\Big(T(\log T)^{k^2-1+\varepsilon}\Big).
\end{equation}
Analogously, the general term turns out to be
\begin{equation}\begin{split}\notag
\int_T^{2T}\Big(\frac{t}{2\pi}\Big)^{-\mathcal S-\mathcal T}\partial_z^n\bigg[e^{\frac{z^2}{4}\mathcal L-k\mu z}\exp\bigg(\frac{z}{2}\sum_{p\in X}&\frac{g_p(\mathcal S,\mathcal T)}{p}\bigg) \bigg]_{z=0}M_{\alpha_{\mathcal S},\beta_{\mathcal T}}(0) dt
+ O_{k,n}\Big(T(\log T)^{k^2-1+\varepsilon}\Big)
\end{split}\end{equation}
and putting this into (\ref{P3.2}), i.e. summing over $j,\mathcal S,\mathcal T$, we get the thesis.
\endproof

\subsection*{Proof of Theorem 1} This proof follows easily from the three propositions we have proved above. If we take $\mu=\mathcal L$ and all the shifts small enough, i.e. $\alpha_i,\beta_j\ll(\log T)^{-1}$ for any $i,j$, then the exponent in the right hand side of the thesis in Proposition \ref{prop3} becomes
\begin{equation}\begin{split}\notag
\frac{z^2}{4}\mathcal L+\frac{z}{2}\sum_{p\in X}\frac{g_p(\mathcal S,\mathcal T)}{p}-zk\mathcal L
&= \frac{z^2}{4}\mathcal L+O\bigg(z\sum_{p\in X}\frac{|p^{-\alpha_1}-1|+\dots+|p^{-\beta_k}-1|}{p}\bigg) \\&
= \frac{z^2}{4}\mathcal L+O_k\bigg(z\sum_{p\in X}\frac{|p^{-\alpha_1}-1|}{p}\bigg)
=\frac{z^2}{4}\mathcal L+O_k(z)
\end{split}\end{equation}
and does not depend on $\mathcal S$ and $\mathcal T$ asymptotically, in fact. Hence we can bring that factor outside and reconstruct the moment of $\zeta_{\alpha,\beta}$, as follows
\begin{equation}\begin{split}\notag
&\int_T^{2T}(\Re P(t)-k\mathcal L)^n\zeta_{\alpha,\beta}(t)dt \\&\hspace{.3cm}
=\partial_z^n\Big[e^{\frac{z^2}{4}\mathcal L+O_k(z)}\Big]_{z=0}
\int_T^{2T}\sum_{0\leq j\leq k}\sum_{\substack{\mathcal S\in\Phi_j \\ \mathcal T\in\Psi_j}}\Big(\frac{t}{2\pi}\Big)^{-\mathcal S-\mathcal T}M_{\alpha_{\mathcal S},\beta_{\mathcal T}}(0)
dt +O_{k,n}\Big(T(\log T)^{k^2-1+\varepsilon}\Big) \\&\hspace{.3cm}
=\bigg( \partial_z^n\Big[e^{\frac{z^2}{4}\mathcal L}\Big]_{z=0}+O_{k,n}\Big((\log\log T)^{\frac{n-1}{2}}\Big)\bigg)
\int_T^{2T}\zeta_{\alpha,\beta}(t)dt + O_{k,n}\Big(T(\log T)^{k^2-1+\varepsilon}\Big).
\end{split}\end{equation}
The thesis follows by analyzing 
\begin{equation}\label{gaussiancoefficient}\partial_z^n \Big[e^{\frac{z^2}{4}\mathcal L}\Big]_{z=0}=\bigg[\sum_{m_1+2m_2=n}\frac{n!}{m_1!m_2!}\bigg(\frac{2z\mathcal L}{4}\bigg)^{m_1}\bigg(\frac{1}{2!}\frac{\mathcal L}{2}\bigg)^{m_2}\bigg]_{z=0}.\end{equation}
If $n$ is odd then the coefficient $\partial_z^n [e^{\frac{z^2}{4}\mathcal L}]_{z=0}$ vanishes, while if $n$ is even then only the term for $m_1=0$ and $m_2=\frac{n}{2}$ survives and gives $\frac{n!}{(n/2)!}(\frac{\mathcal L}{4})^{n/2}=(n-1)!!(\frac{\mathcal L}{2})^{n/2}$, i.e.
\begin{equation}\notag
\frac{1}{\int_T^{2T}\zeta_{\alpha,\beta}(t)dt}\int_T^{2T}(\Re P(t)-k\mathcal L)^n\zeta_{\alpha,\beta}(t)dt=\begin{cases} (1+o_{k,n}(1))(n-1)!!(\frac{\mathcal L}{2})^{n/2} \quad\text{if }n\text{ even} \\o_{k,n}\big((\log\log T)^{n/2}\big) \hspace{1.9cm}\text{if }n\text{ odd}\end{cases}.
\end{equation}
This matches with the Gaussian coefficient then this proves that, in the limit $T\to \infty$, $\Re P(t)$ has Gaussian distribution, with mean $k\log\log T$ and variance $\frac{1}{2}\log\log T$ and then so does $\log|\zeta(1/2+it)|$, in view of Proposition \ref{prop1} and \ref{prop2}. Theorem \ref{caso2k} follows by taking the derivatives with respect to the shifts.

\subsection*{Proof of Theorem 2} To derive Theorem \ref{thmshift}, in Proposition \ref{prop3} we set $\alpha_1=\dots=\alpha_k=i\alpha$ and $\beta_1=\dots=\beta_k=-i\alpha$, with $\alpha\in\mathbb R, |\alpha|<1$ and we take $\mu$ as
\begin{equation}\label{mediashift} \mu_\alpha:=\sum_{p\in X}\frac{\cos(\alpha\log p)}{p} = \widetilde\mu_{\alpha} +O(1)=\begin{cases}\log\log T+O(1) \quad\; \text{if } |\alpha|\log T\leq 1\\-\log|\alpha|+O(1) \quad\; \text{if } |\alpha|\log T>1 \end{cases}\end{equation} by partial summation. Then we get
\begin{equation}\begin{split}\notag
&\int_T^{2T}(\Re P(t)-k\mu_\alpha)^n|\zeta(1/2+i\alpha+it)|^{2k}dt \\&\hspace{3cm}
=\big(1+o_{k,n}(1)\big) \; \partial_z^n\Big[e^{\frac{z^2}{4}\mathcal L}\Big]_{z=0}\int_T^{2T}|\zeta(1/2+i\alpha+it)|^{2k}dt
\end{split}\end{equation}
since the quantity $\frac{z}{2}\sum_{p\in X}\frac{g_p(\mathcal S,\mathcal T)}{p}-zk\mu_\alpha$ vanishes for all $\mathcal S,\mathcal T$ for this choice of the shifts. The thesis follows as in the proof of Theorem \ref{caso2k}.


\section{Proof of Theorem \ref{rmt}} 

In the usual notations we set in the introduction, let us define the moment generating function
\begin{equation}\label{mgf} M_N(s)=\langle |Z|^s\rangle=\sum_{j=0}^\infty\frac{\langle(\log |Z|)^j\rangle}{j!}s^j\end{equation}
(where the mean has to be considered over the group $U(N)$ with respect to the Haar measure) and the cumulants $Q_j=Q_j(N)$ by
\begin{equation}\label{cumulants} \log M_N(s)=\sum_{j=1}^\infty\frac{Q_j}{j!}s^j.\end{equation}
In \cite{RM2000}, among other things, Keating and Snaith studied the cumulants showing that 
\begin{equation}\label{cum}Q_n=\begin{cases}0\hspace{2.9cm}\text{ if }n=1\\ \frac{1}{2}\log N+O(1)\hspace{0.5cm}\text{ if }n=2\\ O(1)\hspace{2.3cm} \text{ if }n\geq 3\end{cases} \end{equation}
 and deduced a central limit theorem proving that the limiting distribution of $\log |Z|$ is Gaussian with mean $0$ and variance $\frac{1}{2}\log N$. Here, for any $k\in\mathbb N$, we study the distribution of random variable $\log |Z|$ with respect to the tilted measure $|Z|^{2k}d_{Haar}$. Before starting with our analysis, we recall that the moments of $|Z|$ are known also for non integer $k$ (see \cite{RM2000} equation (6) and (16)):
\begin{equation}\label{omega}M_N(2k)=\langle |Z|^{2k}\rangle= \exp\bigg(\sum_{j=1}^\infty \frac{(2k)^j}{j!}Q_j\bigg)=\prod_{j=1}^N\frac{\Gamma(j)\Gamma(j+2k)}{\Gamma(j+k)^2}\sim N^{k^2}\frac{G^2(1+k)}{G(1+2k)} \end{equation}
where $k\in\mathbb R$ and $G$ denotes the Barnes $G$-function. We denote $ \mathcal M_{2k}:= \prod_{j=1}^N\frac{\Gamma(j)\Gamma(j+2k)}{\Gamma(j+k)^2}$.\\

Now we are ready to consider the first moment
\begin{equation}\begin{split}\notag
\langle |Z|^{2k}\log|Z|\rangle =\frac{d}{dx}\bigg[\langle|Z|^{2k+x}\rangle\bigg]_{x=0}= \frac{d}{dx}\bigg[\prod_{j=1}^N\frac{\Gamma(j)\Gamma(j+2k+x)}{\Gamma(j+k+\frac{x}{2})^2}\bigg]_{x=0}
\end{split}\end{equation}
by (\ref{mgf}) and (\ref{omega}). We compute the derivative by Leibniz's rule, writing 
$$\prod_{j=1}^N\frac{\Gamma(j)\Gamma(j+2k+x)}{\Gamma(j+k+\frac{x}{2})^2}=\exp\bigg({\sum_{j=1}^N\Big\{\log\Gamma(j)+\log\Gamma(j+2k+x)-2\log\Gamma(j+k+x/2)\Big\}}\bigg)$$
and we get
\begin{equation}\begin{split}\label{aux5.1}
\langle |Z|^{2k}\log|Z|\rangle =\prod_{j=1}^N\frac{\Gamma(j)\Gamma(j+2k)}{\Gamma(j+k)^2}\sum_{j=1}^N\bigg\{\frac{\Gamma'}{\Gamma}(j+2k)-\frac{\Gamma'}{\Gamma}(j+k)\bigg\}.
\end{split}\end{equation}
Moreover an application of Stirling's formula yields 
\begin{equation}\begin{split}\label{aux5.2}
\frac{\Gamma'}{\Gamma}(j+2k)-\frac{\Gamma'}{\Gamma}(j+k)=\frac{k}{j}+O_{k}\bigg(\frac{1}{j^2}\bigg)
\end{split}\end{equation}
hence, by (\ref{aux5.1}) and (\ref{aux5.2}), the weighted mean of the random variable $\log |Z|$ is 
\begin{equation}\begin{split}\notag
\mu_{2k}:=\frac{1}{\mathcal M_{2k}}\langle |Z|^{2k}\log|Z|\rangle =k\log N+O_{k}(1).
\end{split}\end{equation}
Then we study the weighted $n$-th moment of the random variable $\log |Z|$:
\begin{equation}\begin{split}\label{aux5.3}
\langle |Z|^{2k}(\log &|Z|-\mu_{2k})^n\rangle=\sum_{h+j=n}\binom{n}{h}(-\mu_{2k})^j\langle |Z|^{2k}(\log|Z|)^h\rangle \\&
=\sum_{h+j=n}\binom{n}{h}\frac{d^j}{dx^j}\bigg[e^{-x\mu_{2k}}\bigg]_{x=0}\frac{d^h}{dx^h}\bigg[\prod_{j=1}^N\frac{\Gamma(j)\Gamma(j+2k+x)}{\Gamma(j+k+\frac{x}{2})^2}\bigg]_{x=0}\\&
=\frac{d^n}{dx^n}\bigg[\exp\bigg(-x\mu_{2k}+\sum_{j=1}^N\log\bigg(\frac{\Gamma(j)\Gamma(j+2k+x)}{\Gamma(j+k+\frac{x}{2})^2}\bigg)\bigg)\bigg]_{x=0}.
\end{split}\end{equation}
If we denote $f_j(x)=f_{N,k,j}(x):=\log\Gamma(j)+\log\Gamma(j+2k+x)-2\log\Gamma(j+k+\frac{x}{2})$, then we can carry on the computation in (\ref{aux5.3}) by computing the derivative, getting
\begin{equation}\begin{split}\label{aux5.4}
&\frac{d^n}{dx^n}\bigg[\exp\bigg(-x\mu_{2k}+\sum_{j=1}^N\log\bigg(\frac{\Gamma(j)\Gamma(j+2k+x)}{\Gamma(j+k+\frac{x}{2})^2}\bigg)\bigg)\bigg]_{x=0} \\&
=\sum_{m_1,\dots,m_n}\frac{n!}{m_1!\cdots m_n!}e^{\sum_{j=1}^N f_j(0)}\prod_{i=1}^n \bigg(\frac{1}{i!}\frac{d^i}{dx^i}\bigg[-x\mu_{2k}+\sum_{j=1}^Nf_j(x)\bigg]_{x=0}\bigg)^{m_i}\\&
=\mathcal M_{2k}\sum_{m_1,\dots,m_n}\frac{n!}{m_1!\cdots m_n!}\bigg(-\mu_{2k}+\sum_{j=1}^Nf_j'(0)\bigg)^{m_1}\bigg(\frac{1}{2}\sum_{j=1}^Nf_j''(0)\bigg)^{m_2}\prod_{i=3}^n\bigg(\frac{1}{i!}\sum_{j=1}^Nf_j^{(i)}(0)\bigg)^{m_i}
\end{split}\end{equation}
where the sums in (\ref{aux5.4}) are over the $n$-uple $(m_1,\dots,m_n)$ such that $$m_1+2m_2+\cdots+nm_n=n.$$
Using Stirling's approximation formula, one can easily estimate the derivatives of $f_j(x)$ and prove
\begin{equation}\begin{split}\label{aux5.5}
&\sum_{j=1}^Nf_j'(0)=\sum_{j=1}^N\Big\{\frac{k}{j}+O_k(j^{-2})\Big\}=k\log N+O_k(1);
\\&\sum_{j=1}^Nf_j''(0)=\sum_{j=1}^N\Big\{\frac{1}{j+2k}-\frac{1/2}{j+k}+O(j^{-2})\Big\}=\frac{1}{2}\log N+O_k(1);
\\&\sum_{j=1}^Nf_j^{(i)}(0)=\sum_{j=1}^NO(j^{-2})=O(1) \quad\text{for all }i\geq 3.
\end{split}\end{equation}
Putting together (\ref{aux5.3}), (\ref{aux5.4}) and (\ref{aux5.5}) one has
\begin{equation}\begin{split}\notag
&\langle |Z|^{2k}(\log |Z|-\mu_{2k})^n\rangle\\&
=\mathcal M_{2k}\sum_{m_1+2m_2+\cdots+nm_n=n}\frac{n!}{m_1!\cdots m_n!}\bigg(\frac{1}{4}\log N+O_k(1)\bigg)^{m_2}\bigg(O_k(1)\bigg)^{m_1+m_3+\cdots+m_n}
\end{split}\end{equation}
then if $n$ is even the asymptotic is given by $m_2=n/2$ and $m_i=0$ for $i\not =2$, giving
\begin{equation}\begin{split}\notag
\langle |Z|^{2k}(\log |Z|-\mu_{2k})^n\rangle\sim_{k,n} \mathcal M_{2k}\frac{n!}{(n/2)!}\Big(\frac{1}{4}\log N\Big)^{n/2}=\mathcal M_{2k}(n-1)!!\Big(\frac{1}{2}\log N\Big)^{n/2}
\end{split}\end{equation}
while if $n$ is odd the $n$-th moment is surely $o_{k,n}(\mathcal M_{2k}(\log N)^{n/2})$.\\

\textbf{Acknowledgments}.  I am grateful to Sandro Bettin for his support and encouragement as well as for many useful suggestions and to Jon Keating for helpful conversations and for pointing me out interesting connections with the random matrix theory setting, which inspired, among other things, the last section of this paper. I also wish to thank the referee for a very careful reading of the paper and for indicating several inaccuracies.


{\small

\end{document}